\documentclass[12pt,draftclsnofoot,onecolumn]{IEEEtran}

\usepackage{comment,url}
\usepackage[dvips]{graphicx}
\usepackage{caption,subfig}
\usepackage{float}
\graphicspath{{../fig/}}
\DeclareGraphicsExtensions{.eps}
\usepackage{setspace}
\usepackage{amsmath}
\usepackage{amssymb}
\usepackage{cite}
\usepackage{array}
\usepackage{mdwmath}
\usepackage{mdwtab}

\begin{document}

\title{Random Distances Associated with Rhombuses} 

\author{Yanyan Zhuang and Jianping Pan\\
University of Victoria, Victoria, BC, Canada}

\maketitle

\begin{abstract}
Parallelograms are one of the basic building blocks in two-dimensional tiling.
They have important applications in a wide variety of science and engineering
fields, such as wireless communication networks, urban transportation,
operations research, etc. Different from rectangles and squares, the coordinates
of a random point in parallelograms are no longer independent. As a
case study of parallelograms, the explicit probability density functions of the
random Euclidean distances associated with rhombuses are given in this report,
when both endpoints are randomly distributed in 1) the same rhombus, 2) two
parallel rhombuses sharing a side, and 3) two rhombuses having a common
diagonal, respectively. The accuracy of the distance distribution functions is verified
by simulation, and the correctness is validated by a recursion and a
probabilistic sum. The first two statistical moments of the random distances,
and the polynomial fit of the density functions are also given in this report for practical uses. 
\end{abstract}


\section{The Problem}

\begin{figure}
  \centering
  \includegraphics[width=0.5\columnwidth]{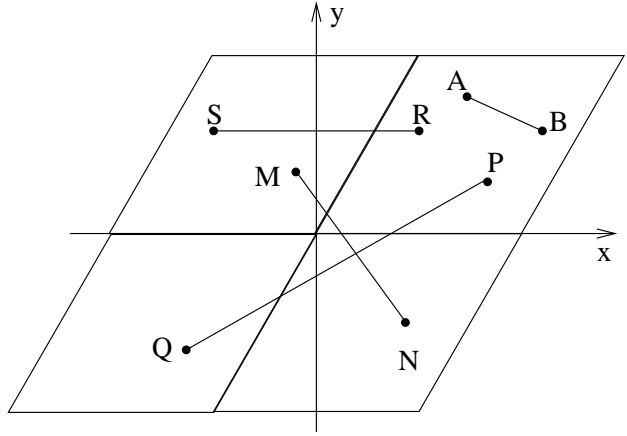}
  \caption{Random Points and Distances Associated with Rhombuses.}
  \label{fig:rhombus}
\end{figure}

Define a ``unit rhombus'' as the rhombus with an acute angle of
$\theta=\frac{\pi}{3}$ and a side length of $1$. Picking two points
uniformly at random from the interior of a unit rhombus, or between two adjacent
unit
rhombuses, the goal is to obtain the probabilistic density function (PDF) of
the random distances between these two endpoints, as illustrated in
Fig.~\ref{fig:rhombus}. 

There are four cases depending on the geometric locations of these two random
endpoints, when
rhombuses are adjacent and similarly oriented, as shown in
Fig.~\ref{fig:rhombus}: i.e., ${\rm AB}$ that are within the same rhombus;
${\rm RS}$ that are inside two parallel rhombuses sharing a side; ${\rm PQ}$
and ${\rm MN}$ that are inside two rhombuses sharing a common diagonal. Here
${\rm PQ}$ and ${\rm MN}$ are two different cases, and in the following, we
refer to them as long (long-diag) and short diagonal (short-diag),
respectively. The next section gives the explicit PDFs for the above four cases.

\section{Distance Distributions Associated with Rhombuses}\label{sec:result}
\subsection{$|AB|$: Distance Distribution within a Rhombus}
{\rm The probability density function of the random Euclidean distances
    between two uniformly distributed points that are both inside the same
    unit rhombus is}
 \begin{equation}\label{eq:fd_r_within}
  f_{D_{\rm I}}(d)=2d\left\{
    \begin{array}{lr}

\left(\frac{4}{3}+\frac{2\pi}{9\sqrt{3}}\right)d^2-\frac{16}{3}d+\frac{2\pi}{
\sqrt{3}} & 0\leq d\leq \frac{\sqrt{3}}{2}\\

\frac{8}{\sqrt{3}}\left(1+\frac{d^2}{3}\right)\sin^{-1}\frac{\sqrt{3}}{2d}
+\left(\frac{4}{3}-\frac{10\pi}{9\sqrt{3}}\right)d^2-\frac{16}{3}d+\frac{10}{3}
\sqrt{4d^2-3}-\frac{2\pi}{\sqrt{3}} & \frac{\sqrt{3}}{2}\leq d\leq 1\\

\frac{4}{\sqrt{3}}\left(1-\frac{d^2}{3}\right)\sin^{-1}\frac{\sqrt{3}}{2d}
-\left(\frac{2}{3}-\frac{2\pi}{9\sqrt{3}}\right)d^2+\sqrt{4d^2-3}-\frac{2\pi}{
3\sqrt{3}}-1 & 1\leq d\leq \sqrt{3} \\

      0 & {\rm otherwise}
    \end{array}
  \right..
\end{equation}

The corresponding cumulative distribution function (CDF) is 
\begin{equation}\label{eq:Fd_r_within}
  F_{D_{\rm I}}(d)=\left\{
    \begin{array}{lr}

\left(\frac{2}{3}+\frac{\pi}{9\sqrt{3}}\right)d^4-\frac{32}{9}d^3+\frac{
2\pi}{\sqrt{3}}d^2 & 0\leq d\leq \frac{\sqrt{3}}{2}\\

\frac{4}{\sqrt{3}}\left(2d^2+\frac{d^4}{3}\right)\sin^{-1}\frac{\sqrt{3}}{
2d}+\left(\frac{2}{3}-\frac{5\pi}{9\sqrt{3}}\right)d^4-\frac{32
}{9}d^3-\frac{2\pi}{\sqrt{3}}d^2 \\
~~~~+\frac{14d^2+3}{6}\sqrt{4d^2-3} & \frac{\sqrt{3}}{2}\leq d\leq 1\\

\frac{2}{\sqrt{3}}\left(2d^2-\frac{d^4}{3}\right)\sin^{-1}\frac{\sqrt{3}}{
2d}+\left(\frac{\pi}{9\sqrt{3}}-\frac{1}{3}\right)d^4-\left(\frac{2\pi}{3\sqrt{3
}}+1\right)d^2 \\
~~~~+\frac{22d^2+15}{36}\sqrt{4d^2-3}+\frac{1}{4} & 1\leq d\leq \sqrt{3} \\

0 & {\rm otherwise}
    \end{array}
 \right..
\end{equation}

\subsection{$|RS|$: Distance Distribution between Two Parallel Rhombuses Sharing
a Side}
{\rm The probability density function of the random distances
    between two uniformly distributed points, one in each of the two adjacent 
    unit rhombuses that are parallel to each other, is}
 \begin{equation}\label{eq:fd_r_para}
  f_{D_{\rm P}}(d)=2d\left\{
    \begin{array}{lr}

\frac{4}{3}d-\left(\frac{2}{3}+\frac{\pi}{9\sqrt{3}}\right)d^2 & 0\leq d\leq
\frac{\sqrt{3}}{2}\\

-\frac{2}{\sqrt{3}}(d^2+2)\sin^{-1}\frac{\sqrt{3}}{2d}+\left(\frac{8\pi}
{9\sqrt{3}}-\frac{2}{3}\right)d^2+\frac{4}{3}d-\frac{11}{6}\sqrt{4d^2-3} \\
~~~~+\frac{2\pi}{\sqrt{3}} & \frac{\sqrt{3}}{2}\leq d\leq 1\\

\frac{4d^2}{3\sqrt{3}}\sin^{-1}\frac{\sqrt{3}}{2d}+\left(\frac{2}{3}
-\frac{2\pi}{9\sqrt{3}}\right)d^2-\frac{8}{3}d+\frac{\sqrt{4d^2-3}}{3}+\frac{
2\pi}{3\sqrt{3}}+\frac{1}{2} & 1\leq d\leq \sqrt{3} \\

\left(\frac{2}{\sqrt{3}}-\frac{d^2}{3\sqrt{3}}\right)\sin^{-1}\frac{\sqrt
{3}}{2d}+\left(\frac{2}{\sqrt{3}}+\frac{d^2}{3\sqrt{3}}\right)\sin^{-1}\frac{
\sqrt{3}}{d}+\left(\frac{1}{3}-\frac{\pi}{9\sqrt{3}}\right)d^2 \\
~~~~-\frac{8}{3}d+\frac{7}{12}\sqrt{4d^2-3}+\sqrt{d^2-3}+\frac{3}{4}
-\frac{2\pi}{3\sqrt{3}} & \sqrt{3}\leq d\leq 2 \\

\left(\frac{2}{\sqrt{3}}-\frac{d^2}{3\sqrt{3}}\right)\left(\sin^{-1}\frac{
\sqrt{3}}{2d}+\sin^{-1}\frac{\sqrt{3}}{d}\right)+\left(\frac{\pi}{9\sqrt{3}}
-\frac{1}{3 }\right)d^2 \\
~~~~+\frac{7}{12}\sqrt{4d^2-3}+\frac{\sqrt{d^2-3}}{3}-\frac{2\pi}{3\sqrt{3}}
-\frac{5}{4} & 2 \leq d \leq \sqrt{7} \\

      0 & {\rm otherwise}
    \end{array}
 \right..
\end{equation}

The corresponding CDF is 
\begin{equation}\label{eq:Fd_r_para}
F_{D_{\rm P}}(d)=\left\{
    \begin{array}{lr} 

\frac{8}{9}d^3-\left(\frac{1}{3}+\frac{\pi}{18\sqrt{3}} \right)d^4 & 0\leq
d\leq \frac{\sqrt{3}}{2}\\

-\frac{1}{\sqrt{3}}\left(4d^2+d^4\right)\sin^{-1}\frac{\sqrt{3}}{
2d}+\left(\frac{4\pi}{9\sqrt{3}}-\frac{1}{3}\right)d^4+\frac{8}{9}d^3
+\frac{2\pi}{\sqrt{3}}d^2 \\
~~~~-\frac{94d^2+15}{72}\sqrt{4d^2-3} & \frac{\sqrt{3}}{2}\leq d\leq 1\\

\frac{2d^4}{3\sqrt{3}}\sin^{-1}\frac{\sqrt{3}}{2d}+\left(\frac{1}{3}
-\frac{\pi}{9\sqrt{3}}\right)d^4-\frac{16}{9}d^3+\left(\frac{2\pi}{3\sqrt{3}
}+\frac{1}{2}\right)d^2 \\
~~~~+\frac{10d^2-3}{36}\sqrt{4d^2-3}-\frac{5}{24} & 1\leq d\leq \sqrt{3} \\

\left(\frac{2d^2}{\sqrt{3}}-\frac{d^4}{6\sqrt{3}}\right)\sin^{-1}\frac{\sqrt{3}
}{2d}+\left(\frac{2d^2}{\sqrt{3}}+\frac{d^4}{6\sqrt{
3}}\right)\sin^{-1}\frac{\sqrt{3}}{d}+\left(\frac{1}{6}-\frac{\pi}{18\sqrt{3}}
\right)d^4 \\
~~~~-\frac{16}{9}d^3+\left(\frac{3}{4}-\frac{2\pi}{3\sqrt{3}}\right)d^2+\frac{
6d^2+3}{16}\sqrt{4d^2-3}+\frac{13d^2+6}{18}\sqrt{d^2-3}-\frac{55}{48} &
\sqrt{3}\leq d\leq 2 \\

\left(\frac{2d^2}{\sqrt{3}}-\frac{d^4}{6\sqrt{3}}\right)\left(\sin^{-1}\frac{
\sqrt{3}}{2d}+\sin^{-1}\frac{\sqrt{3}}{d}\right)+\left(\frac{\pi}{18\sqrt
{3}}-\frac{1}{6}\right)d^4 \\
~~~~-\left(\frac{2\pi}{3\sqrt{3}}+\frac{5}{4}\right)d^2+\frac{6d^2+3}{16}\sqrt{
4d^2-3}+\frac{d^2+6}{6}\sqrt{d^2-3}-\frac{23}{48} & 2 \leq d \leq \sqrt{7} \\

     0 & {\rm otherwise}
    \end{array}
 \right..
\end{equation}

\subsection{$|PQ|$: Distance Distribution between Two Long-Diag Adjacent
Rhombuses}
{\rm The probability density function of the random distances between two
    uniformly distributed points, one in each of the two adjacent unit rhombuses
    that have a common long diagonal, is}
 \begin{equation}\label{eq:fd_r_diag1}
  f_{D_{\rm LD}}(d)=2d\left\{
    \begin{array}{lr}

\left(\frac{1}{3}-\frac{\pi}{9\sqrt{3}}\right)d^2 & 0\leq d\leq 1\\

-\frac{4d^2}{3\sqrt{3}}\sin^{-1}\frac{\sqrt{3}}{2d}+\left(\frac{\pi}{3\sqrt{3}}
-1\right)d^2+\frac{8}{3}d-\frac{\sqrt{4d^2-3}}{3}-1 & 1\leq d\leq \sqrt{3}\\

\frac{4}{\sqrt{3}}\left(\frac{d^2}{3}-2\right)\sin^{-1}\frac{\sqrt{3
}}{2d}+\left(\frac{1}{3}-\frac{\pi}{9\sqrt{3}}\right)d^2+\frac{8}{3}d-\frac{7}{3
}\sqrt{4d^2-3}+\frac{4\pi}{3\sqrt{3}}+1 & \sqrt{3}\leq d\leq 2\\

\frac{4}{\sqrt{3}}\left(\frac{d^2}{3}-2\right)\sin^{-1}\frac{\sqrt{3
}}{2d}+\frac{2d^2}{3\sqrt{3}}\sin^{-1}\frac{\sqrt{3}}{d}+\left(1-\frac{\pi}{
3\sqrt{3}}\right)d^2 \\
~~~~-\frac{7}{3}\sqrt{4d^2-3}+\frac{2}{3}\sqrt{d^2-3}+\frac{4\pi}{3\sqrt{3}}+3 &
2\leq d\leq \sqrt{7} \\

\frac{2}{\sqrt{3}}\left(4-\frac{d^2}{3}\right)\sin^{-1}\frac{\sqrt{3}}{d}
+\left(\frac{\pi}{9\sqrt{3}}-\frac{1}{3}\right)d^2+2\sqrt{d^2-3}-\frac{4\pi}{
3\sqrt{3}}-2 & \sqrt{7} \leq d \leq 2\sqrt{3} \\

      0 & {\rm otherwise}
    \end{array}
  \right..
\end{equation}

The corresponding CDF is 
\begin{equation}\label{eq:Fd_r_diag1}
  F_{D_{\rm LD}}(d)=\left\{
    \begin{array}{lr}

\left(\frac{1}{6}-\frac{\pi}{18\sqrt{3}}\right)d^4 & 0\leq d\leq 1\\

-\frac{2d^4}{3\sqrt{3}}\sin^{-1}\frac{\sqrt{3}}{2d}+\left(\frac{\pi}{6\sqrt{3}
}-\frac{1}{2}\right)d^4+\frac{16}{9}d^3-d^2-\frac{10d^2-3}{36}\sqrt{4d^2-3}
+\frac{1}{12} & 1\leq d\leq \sqrt{3}\\

\frac{2}{\sqrt{3}}\left(\frac{d^4}{3}-4d^2\right)\sin^{-1}\frac{\sqrt{3
}}{2d}+\left(\frac{1}{6}-\frac{\pi}{18\sqrt{3}}\right)d^4+\frac{16}{9}d^3
+\left(1+\frac{4\pi}{3\sqrt{3}}\right)d^2 \\
~~~~-\frac{6d^2+3}{4}\sqrt{4d^2-3}+\frac{19}{12} & \sqrt{3}\leq d\leq 2\\

\frac{2}{\sqrt{3}}\left(\frac{d^4}{3}-4d^2\right)\sin^{-1}\frac{\sqrt{3
}}{2d}+\frac{d^4}{3\sqrt{3}}\sin^{-1}\frac{\sqrt{3}}{d}+\left(\frac{1}{2}-\frac{
\pi}{6\sqrt{3}}\right)d^4 \\
~~~~+\left(3+\frac{4\pi}{3\sqrt{3}}\right)d^2-\frac{6d^2+3}{4}\sqrt{4d^2-3}
+\frac{5d^2-6}{9}\sqrt{d^2-3}+\frac{11}{12} & 2\leq d\leq \sqrt{7} \\

\frac{1}{\sqrt{3}}\left(8d^2-\frac{d^4}{3}\right)\sin^{-1}\frac{\sqrt{3}
}{d}+\left(\frac{\pi}{18\sqrt{3}}-\frac{1}{6}\right)d^4-\left(\frac{4\pi}{3\sqrt
{3}}+2\right)d^2 \\
~~~~+\frac{11d^2+30}{9}\sqrt{d^2-3}-5 & \sqrt{7} \leq d \leq 2\sqrt{3} \\

      0 & {\rm otherwise}
    \end{array}
  \right..
\end{equation}

\subsection{$|MN|$: Distance Distribution between Two Short-Diag Adjacent
Rhombuses}
{\rm The probability density function of the random distances between two
    uniformly distributed points, one in each of the two adjacent unit rhombuses
    that have a common short diagonal, is}
 \begin{equation}\label{eq:fd_r_diag2}
  f_{D_{\rm SD}}(d)=2d\left\{
    \begin{array}{lr}

\left(\frac{1}{3}+\frac{2\pi}{9\sqrt{3}}\right)d^2 & 0\leq d\leq
\frac{\sqrt{3}}{2}\\

\frac{8d^2}{3\sqrt{3}}\sin^{-1}\frac{\sqrt{3}}{2d}+\left(\frac{1}{3}-\frac{10\pi
}{9\sqrt{3}}\right)d^2+\frac{2}{3}\sqrt{4d^2-3} & \frac{\sqrt{3}}{2} \leq d \leq
1\\

-\frac{4}{\sqrt{3}}\left(\frac{d^2}{3}+2\right)\sin^{-1}\frac{\sqrt{3}}{
2d}+\left(\frac{1}{3}+\frac{2\pi}{9\sqrt{3}}\right)d^2+\frac{8}{3}d-3\sqrt{
4d^2-3}+\frac{8\pi}{3\sqrt{3}}+1 & 1\leq d\leq \sqrt{3} \\

\frac{8}{\sqrt{3}}\sin^{-1}\frac{\sqrt{3}}{d}-d^2+\frac{8}{3}d+\frac{8}{3}\sqrt{
d^2-3}-\frac{8\pi}{3\sqrt{3}}-4 & \sqrt{3}\leq d\leq 2 \\

      0 & {\rm otherwise}
    \end{array}
  \right..
\end{equation}

The corresponding CDF is 
 \begin{equation}\label{eq:Fd_r_diag2}
  F_{D_{\rm SD}}(d)=\left\{
    \begin{array}{lr}

\left(\frac{1}{6}+\frac{\pi}{9\sqrt{3}}\right)d^4 & 0\leq d\leq
\frac{\sqrt{3}}{2}\\

\frac{4d^4}{3\sqrt{3}}\sin^{-1}\frac{\sqrt{3}}{2d}+\left(\frac{1}{6}-\frac{5\pi}
{9\sqrt{3}}\right)d^4+\frac{10d^2-3}{18}\sqrt{4d^2-3} & \frac{\sqrt{3}}{2} \leq
d \leq 1\\

-\frac{2}{\sqrt{3}}\left(\frac{d^4}{3}+4d^2\right)\sin^{-1}\frac{\sqrt{3}}{2d}
+\left(\frac{1}{6}+\frac{\pi}{9\sqrt{3}}\right)d^4+\frac{16}{9}d^3+\left(\frac{
8\pi}{3\sqrt{3}}+1\right)d^2\\
~~~~-\frac{74d^2+21}{36}\sqrt{4d^2-3}+\frac{1}{4} & 1\leq d\leq \sqrt{3} \\

\frac{8d^2}{\sqrt{3}}\sin^{-1}\frac{\sqrt{3}}{d}-\frac{d^4}{2}+\frac{16}{9}
d^3-\left(4+\frac{8\pi}{3\sqrt{3}}\right)d^2+\frac{16d^2+24}{9}\sqrt{d^2-3}+1 &
\sqrt{3}\leq d\leq 2 \\

      0 & {\rm otherwise}
    \end{array}
  \right..
\end{equation}

Note that although unit rhombuses are assumed throughout 
(\ref{eq:fd_r_within})--(\ref{eq:Fd_r_diag2}), the distance distribution
functions can be
easily scaled by a nonzero scalar, for rhombuses of arbitrary side length. For
example, let
the side length of such rhombuses be $s>0$, then
\begin{equation}
 F_{sD}(d)=P(sD\leq d)=P(D\leq \frac{d}{s})=F_D(\frac{d}{s}). \nonumber
\end{equation}
Therefore,
\begin{equation}\label{eq:scale}
 f_{sD}(d)=F'_D(\frac{d}{s})=\frac{1}{s}f_D(\frac{d}{s}).
\end{equation}

\section{Verification and Validation}

\begin{figure}
  \centering
  \includegraphics[width=0.55\columnwidth]{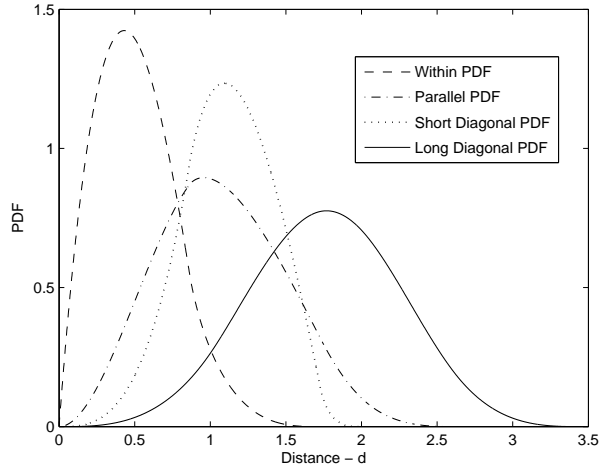}
  \caption{Distributions of Random Distances Associated with
Rhombuses.}
  \label{fig:rhombus_pdf}
\end{figure}

\begin{figure}
  \centering
  \includegraphics[width=0.55\columnwidth]{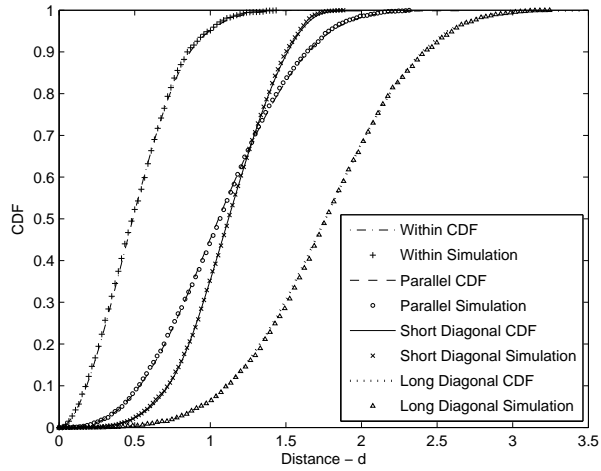}
  \caption{Distribution and Simulation Results For Random Distances Associated with
Rhombuses.}
  \label{fig:rhombus_cdf}
\end{figure}

\subsection{Verification by Simulation}

Figure~\ref{fig:rhombus_pdf} plots the probability density functions, as given in (\ref{eq:fd_r_within}), 
(\ref{eq:fd_r_para}), (\ref{eq:fd_r_diag1}) and (\ref{eq:fd_r_diag2}), 
respectively, of the
four random distance cases shown in Fig.~\ref{fig:rhombus}. Figure~\ref{fig:rhombus_cdf} shows a comparison between the
cumulative distribution functions (CDFs) of the random distances, and the
simulation results by generating $1,000$ pairs of random points with the
corresponding geometric locations as illustrated in
Fig.~\ref{fig:rhombus}. Figure~\ref{fig:rhombus_cdf} demonstrates that our
distance distribution functions are very accurate when compared with the simulation results.

\subsection{Validation by Recursion}

When looking at Fig.~\ref{fig:rhombus}, we find that the four adjacent rhombuses together
resemble a larger rhombus, with a side length of $2$. According to~(\ref{eq:scale}),
the distance distribution in the large rhombus is
$f_{2D}(d)=\frac{1}{2}f_{D_{\rm I}}(\frac{d}{2})$. On the other hand, if we look
at the two random endpoints of a given link inside the large rhombus, they will fall into
one of the four individual cases: both endpoints are inside the same small
rhombus, with probability $\frac{1}{4}$; the two endpoints fall into two
parallel rhombuses, with probability $\frac{1}{2}$; and the two endpoints fall
into two diagonal rhombuses (either long or short-diag), both with probability
$\frac{1}{8}$. Thus the distance density function for the large rhombus can be given by a probabilistic sum,
$f_{2D}(d)=\frac{1}{4}f_{D_{\rm I}}(d)+\frac{1}{2}f_{D_{\rm
P}}(d)+\frac{1}{8}f_{D_{\rm LD}}(d)+\frac{1} {8}f_{D_{\rm SD}}(d)$, where
$f_{D_{\rm I}}(d)$, $f_{D_{\rm P}}(d)$, $f_{D_{\rm LD}}(d)$ and $f_{D_{\rm
SD}}(d)$ are given in (\ref{eq:fd_r_within}), (\ref{eq:fd_r_para}),
(\ref{eq:fd_r_diag1}) and (\ref{eq:fd_r_diag2}), respectively. To confirm that
the above two definitions of $f_{2D}(d)$ are equivalent, i.e., 
$\frac{1}{4}f_{D_{\rm I}}(d)+\frac{1}{2}f_{D_{\rm P}}(d)+\frac{1}{8}f_{D_{\rm
LD}}(d)+\frac{1} {8}f_{D_{\rm SD}}(d)$ is equal to $\frac{1}{2}f_{D_{\rm
I}}(\frac{d}{2})$, we verify them mathematically as follows.

\subsubsection{$0\leq d \leq \frac{\sqrt{3}}{2}$}

\begin{eqnarray}
\frac{1}{4}f_{D_{\rm
I}}(d)=\frac{d}{2}\left[\left(\frac{4}{3}+\frac{2\pi}{9\sqrt{3}}
\right)d^2-\frac{16}{3}d+\frac{2\pi}{\sqrt{3}}\right], \quad
\frac{1}{2}f_{D_{\rm
P}}(d)=d\left[\frac{4}{3}d-\left(\frac{2}{3}+\frac{\pi}{9\sqrt{3}}
\right)d^2\right], \nonumber 
\end{eqnarray}
\begin{eqnarray}
\frac{1}{8}f_{D_{\rm
LD}}(d)=\frac{d}{4}\left(\frac{1}{3}-\frac{\pi}{9\sqrt{3}}\right)d^2, \quad
\frac{1}{8}f_{D_{\rm SD}}(d)=\frac{d}{4}\left(\frac{1}{3}+\frac{2\pi}{9\sqrt{3}}
\right)d^2. \nonumber
\end{eqnarray}

Thus,
 \begin{eqnarray}
f_{2D}(d)&=&\frac{1}{4}f_{D_{\rm I}}(d)+\frac{1}{2}f_{D_{\rm
P}}(d)+\frac{1}{8}f_{D_{\rm LD}}(d)+\frac{1} {8}f_{D_{\rm SD}}(d) \nonumber\\
&=&\frac{d}{2}\left[\left(\frac{4}{3}+\frac{2\pi}{9\sqrt{3}}\right)\left(\frac{d
}{2}\right)^2-\frac{16}{3}\left(\frac{d}{2}\right)+\frac{2\pi}{\sqrt{3}}\right]
=\frac{1}{2}f_{D_{\rm I}}(\frac{d}{2}). \nonumber
 \end{eqnarray}

\subsubsection{$\frac{\sqrt{3}}{2} \leq d \leq 1$}

\begin{eqnarray}
\frac{1}{4}f_{D_{\rm
I}}(d)=\frac{d}{2}\left[\frac{8}{\sqrt{3}}\left(1+\frac{d^2}{3}\right)
\sin^{-1}\frac{\sqrt{3}}{2d}+\left(\frac{4}{3}-\frac{10\pi}{9\sqrt{3}}
\right)d^2-\frac{16}{3}
d+\frac{10}{3}\sqrt{4d^2-3}-\frac{2\pi}{\sqrt{3}}\right], \nonumber 
\end{eqnarray}
\begin{eqnarray}
 \frac{1}{2}f_{D_{\rm P}}(d)=d\left[-\frac{2}{\sqrt{3}}(d^2+2)\sin^{-1}
\frac{\sqrt{3}}{2d}+\left(\frac{8\pi}{9\sqrt{3}}-\frac{2}{3}\right)d^2+
\frac{4}{3}d-\frac{11}{6}\sqrt{4d^2-3}+\frac{2\pi}{\sqrt{3}}\right], \nonumber 
\end{eqnarray}
\begin{eqnarray}
\frac{1}{8}f_{D_{\rm
LD}}(d)=\frac{d}{4}\left(\frac{1}{3}-\frac{\pi}{9\sqrt{3}}\right)d^2, \quad
\frac{1}{8}f_{D_{\rm SD}}(d)=\frac{d}{4}\left[\frac{8d^2}{3\sqrt{3}}\sin^{-1}
\frac{\sqrt{3}}{2d}+\left(\frac{1}{3}-\frac{10\pi}{9\sqrt{3}}\right)d^2+\frac{2}
{3} \sqrt{4d^2-3}\right]. \nonumber
\end{eqnarray}

Thus, 
 \begin{eqnarray}
f_{2D}(d)&=&\frac{1}{4}f_{D_{\rm I}}(d)+\frac{1}{2}f_{D_{\rm
P}}(d)+\frac{1}{8}f_{D_{\rm LD}}(d)+\frac{1} {8}f_{D_{\rm
SD}}(d) \nonumber\\
&=&\frac{d}{2}\left[\left(\frac{4}{3}+\frac{2\pi}{9\sqrt{3}}
\right)\left(\frac{d}{2}\right)^2-\frac{16}{3}\left(\frac{d}{2}\right)+\frac{
2\pi}{\sqrt{3}}\right]=\frac{1}{2}f_{D_{\rm I}}(\frac{d}{2}). \nonumber
 \end{eqnarray}

\subsubsection{$1 \leq d \leq \sqrt{3}$}

\begin{eqnarray}
\frac{1}{4}f_{D_{\rm I}}(d)=\frac{d}{2}\left[\frac{4}{\sqrt{3}}\left(1-
\frac{d^2}{3}\right)\sin^{-1}\frac{\sqrt{3}}{2d}-\left(\frac{2}{3}
-\frac{2\pi}{9\sqrt{3}}\right)d^2+\sqrt{4d^2-3}-\frac{2\pi}{3\sqrt{3}}
-1\right], \nonumber 
\end{eqnarray}
\begin{eqnarray}
 \frac{1}{2}f_{D_{\rm P}}(d)=d\left[\frac{4d^2}{3\sqrt{3}}\sin^{-1}
\frac{\sqrt{3}}{2d}+\left(\frac{2}{3}-\frac{2\pi}{9\sqrt{3}}\right)d^2
-\frac{8}{3}d+\frac{\sqrt{4d^2-3}}{3}+\frac{2\pi}{3\sqrt{3}}
+\frac{1}{2}\right], \nonumber 
\end{eqnarray}
\begin{eqnarray}
\frac{1}{8}f_{D_{\rm LD}}(d)=\frac{d}{4}\left[-\frac{4d^2}{3\sqrt{3}}
\sin^{-1}\frac{\sqrt{3}}{2d}+\left(\frac{\pi}{3\sqrt{3}}-1\right)d^2+
\frac{8}{3}d-\frac{\sqrt{4d^2-3}}{3}-1\right], \nonumber
\end{eqnarray}
\begin{eqnarray}
\frac{1}{8}f_{D_{\rm SD}}(d)=\frac{d}{4}\left[-\left(\frac{4d^2}{3\sqrt{3}}
+\frac{8}{\sqrt{3}}\right)\sin^{-1}\frac{\sqrt{3}}{2d}+\left(\frac{1}{3}
+\frac{2\pi}{9\sqrt{3}}\right)d^2+\frac{8}{3}d-3\sqrt{4d^2-3}
+\frac{8\pi}{3\sqrt{3}}+1\right]. \nonumber
\end{eqnarray}

Thus, 
 \begin{eqnarray}
f_{2D}(d)&=&\frac{1}{4}f_{D_{\rm I}}(d)+\frac{1}{2}f_{D_{\rm
P}}(d)+\frac{1}{8}f_{D_{\rm LD}}(d)+\frac{1} {8}f_{D_{\rm
SD}}(d) \nonumber\\
&=&\frac{d}{2}\left[\left(\frac{4}{3}+\frac{2\pi}{9\sqrt{3}}
\right)\left(\frac{d}{2}\right)^2-\frac{16}{3}\left(\frac{d}{2}\right)+\frac{
2\pi}{\sqrt{3}}\right]=\frac{1}{2}f_{D_{\rm I}}(\frac{d}{2}). \nonumber
 \end{eqnarray}

\subsubsection{$\sqrt{3} \leq d \leq 2$}

$f_{D_{\rm I}}(d)=0$, and
\begin{eqnarray}
 \frac{1}{2}f_{D_{\rm
P}}(d)&=&d\left[\left(\frac{2}{\sqrt{3}}-\frac{d^2}{3\sqrt{3}}\right)
\sin^{-1}\frac{\sqrt{3}}{2d}+\left(\frac{2}{\sqrt{3}}+\frac{d^2}{3\sqrt{3}}
\right)\sin^{
-1}\frac{\sqrt{3}}{d}+\left(\frac{1}{3}-\frac{\pi}{9\sqrt{3}}\right)d^2-\frac{8}
{3}d\right. \nonumber\\
&&\left.+\sqrt{d^2-3}+\frac{7}{12}\sqrt{4d^2-3}+\frac{3}{4}-\frac{2\pi}{3\sqrt{3
}}\right], \nonumber 
\end{eqnarray}
\begin{eqnarray}
\frac{1}{8}f_{D_{\rm
LD}}(d)=\frac{d}{4}\left[\frac{4}{\sqrt{3}}\left(\frac{d^2}{3}-2\right)\sin^{
-1}\frac{\sqrt{3}}{2d}+\left(\frac{1}{3}-\frac{\pi}{9\sqrt{3}}\right)d^2-\frac{7
}{3}\sqrt{4d^2-3}+\frac{8}{3}d+\frac{4\pi}{3\sqrt{3}}+1\right], \nonumber
\end{eqnarray}
\begin{eqnarray}
\frac{1}{8}f_{D_{\rm SD}}(d)=\frac{d}{4}\left[\frac{8}{
\sqrt{3}}\sin^{-1}\frac{\sqrt{3}}{d}-d^2+\frac{8}{3}d+\frac{8}{3}\sqrt{
d^2-3}-\frac{8\pi}{3\sqrt{3}}-4\right]. \nonumber
\end{eqnarray}

Thus,
 \begin{eqnarray}
f_{2D}(d)&=&\frac{1}{2}f_{D_{\rm P}}(d)+\frac{1}{8}f_{D_{\rm LD}}(d)+\frac{1}
{8}f_{D_{\rm SD}}(d) \nonumber\\
&=&\frac{d}{2}\left[\frac{8}{\sqrt{3}}\left(1+\frac{(d/2)^2}{3}\right)\sin^{-1}
\frac{\sqrt{3}}{2(d/2)}+\left(\frac{4}{3}-\frac{10\pi}{9\sqrt{3}}
\right)\left(\frac{d}{2}\right)^2-\frac{16}{3}\left(\frac{d}{2}
\right)\right. \nonumber\\
&&\left.+\frac{10}{3}\sqrt{4\left(\frac{d}{2}\right)^2-3}-\frac{2\pi}{\sqrt{3}}
\right]=\frac{1}{2}f_{D_{\rm I}}(\frac{d}{2}). \nonumber
 \end{eqnarray}

\subsubsection{$2 \leq d \leq \sqrt{7}$}

$f_{D_{\rm I}}(d)=f_{D_{\rm SD}}(d)=0$, and
\begin{eqnarray}
 \frac{1}{2}f_{D_{\rm
P}}(d)&=&d\left[\left(\frac{2}{\sqrt{3}}-\frac{d^2}{3\sqrt{3}}\right)\left(\sin^
{-1}\frac{\sqrt{3}}{2d}+\sin^{-1}\frac{\sqrt{3}}{d}\right)+\left(\frac{\pi}{
9\sqrt{3}}-\frac{1}{3}\right)d^2+\frac{7}{12}\sqrt{
4d^2-3}\right.\nonumber\\ 
&&\left.+\frac{\sqrt{d^2-3}}{3}-\frac{2\pi}{3\sqrt{3}}-\frac{5}{4}\right],
\nonumber 
\end{eqnarray}
\begin{eqnarray}
\frac{1}{8}f_{D_{\rm
LD}}(d)&=&\frac{d}{4}\left[\frac{4}{\sqrt{3}}\left(\frac{d^2}{3}-2\right)\sin^{
-1} \frac{\sqrt{3}}{2d}+\frac{2d^2}{3\sqrt{3}}\sin^{-1}\frac{\sqrt{3}}{d}
+\left(1-\frac{\pi}{3\sqrt{3}}\right)d^2-\frac{7}{3}\sqrt{4d^2-3}
\right.\nonumber\\ 
&&\left.+\frac{2}{3}\sqrt{d^2-3}+\frac{4\pi}{3\sqrt{3}} +3\right]. \nonumber
\end{eqnarray}

Thus,
 \begin{eqnarray}
f_{2D}(d)&=&\frac{1}{2}f_{D_{\rm P}}(d)+\frac{1}{8}f_{D_{\rm
LD}}(d)\nonumber \\
&=&\frac{d}{2}\left[\frac{4}{\sqrt{3}}\left(1-\frac{(d/2)^2}{3}\right)\sin^{-1}
\frac{\sqrt{3}}{d}-\left(\frac{2}{3}-\frac{2\pi}{9\sqrt{3}}\right)\left(\frac{d}
{2}\right)^2+\sqrt{4\left(\frac{d}{2}\right)^2-3}-\frac{2\pi}{3\sqrt{3}}-1\right
]\nonumber \\
&=&\frac{1}{2}f_{D_{\rm I}}(\frac{d}{2}). \nonumber
 \end{eqnarray}

\subsubsection{$\sqrt{7} \leq d \leq 2\sqrt{3}$}

$f_{D_{\rm I}}(d)=f_{D_{\rm P}}(d)=f_{D_{\rm SD}}(d)=0$, and
 \begin{eqnarray}
f_{2D}(d)&=&\frac{1}{8}f_{D_{\rm
LD}}(d)=\frac{d}{4}\left[\frac{2}{\sqrt{3}}\left(4-\frac{d^2}{3}\right)\sin^{-1}
\frac{\sqrt{3}}{d}+\left(\frac{\pi}{9\sqrt{3}}-\frac{1}{3}\right)d^2+2\sqrt{
d^2-3}-\frac{4\pi}{3\sqrt{3}}-2 \right] \nonumber\\
&=&\frac{d}{2}\left[\frac{4}{\sqrt{3}}\left(1-\frac{(d/2)^2}{3}\right)\sin^{-1}
\frac{\sqrt{3}}{d}-\left(\frac{2}{3}-\frac{2\pi}{9\sqrt{3}}\right)\left(\frac{d}
{2}\right)^2+\sqrt{4\left(\frac{d}{2}\right)^2-3}-\frac{2\pi}{3\sqrt{3}}-1\right
]\nonumber \\
&=&\frac{1}{2}f_{D_{\rm I}}(\frac{d}{2}). \nonumber
 \end{eqnarray}

In summary, we have $f_{2D}(d)=\frac{1}{2}f_{D_{\rm I}}(\frac{d}{2})$ by
recursion, and the probabilistic sum $\frac{1}{4}f_{D_{\rm
I}}(d)+\frac{1}{2}f_{D_{\rm P}}(d)+ \frac{1}{8}f_{D_{\rm LD}}(d)+\frac{1}
{8}f_{D_{\rm SD}}(d)$ is equal to $\frac{1}{2}f_{D_{\rm I}}(\frac{d}{2})$ in all
the cases discussed above. The results are a strong validation of the
correctness of the distance distributions that we have derived.

\section{Practical Results}

\subsection{Statistical Moments of Random Distances}

The distance distribution functions given in Section~\ref{sec:result} can conveniently
lead to all the statistical moments of the random distances associated with rhombuses. Given $f_{D_{\rm
I}}(d)$ in (\ref{eq:fd_r_within}), for example, the first moment (mean) of $d$, i.e., the average distance within
a single rhombus, is 
\begin{eqnarray}
 M_{D_{\rm I}}^{(1)}=\int_0^{\sqrt{3}}xf_{D_{\rm
I}}(x)dx=\frac{\sqrt{3}}{8}+\frac{3}{40}+\frac{1}{80}\left[7\ln
\left(2\sqrt{3}+3\right)-6\ln \left(2\sqrt{3}-3\right)\right]\approx
0.5123783359, \nonumber
\end{eqnarray}
and the second raw moment is
\begin{eqnarray}
 M_{D_{\rm I}}^{(2)}=\int_0^{\sqrt{3}}x^2f_{D_{\rm I}}(x)dx=\frac{1}{3},
\nonumber
\end{eqnarray}
from which the variance (the second central moment) can be derived as
\begin{eqnarray}
 Var_{D_{\rm I}}=M_{D_{\rm I}}^{(2)}-\left[M_{D_{\rm I}}^{(1)}\right]^2\approx
0.0708017742. \nonumber
\end{eqnarray}

When the side length of a rhombus is scaled by $s$, the corresponding first two statistical moments given above then become
\begin{equation}
M_{D_{\rm I}}^{(1)}=0.5123783359s,~~\mbox{}~~M_{D_{\rm I}}^{(2)}=\frac{s}{3}
~~\mbox{ and }~~Var_{D_{\rm I}}=0.0708017742s^2.
\end{equation}

\begin{table}
  \caption{Moments and Variance---Numerical vs Simulation Results}
  \centering
  \begin{tabular}{|c||c|c|c|c|}
    \hline
    Endpoint Geometry & PDF/Sim & $M_{D}^{(1)}$ & $M_{D}^{(2)}$ & $Var_{D}$ \\ \hline \hline
    Within a & $f_{D_{\rm I}}(d)$ & $0.5123783359s$ & $0.3333333333s$ &
$0.0708017742s^2$ \\ 
    \cline{2-5}
    Single Rhombus & Sim & $0.5137344650s$ & $0.3356749448s$ &
$0.0717518443s^2$ \\ \hline
    Between two & $f_{D_{\rm P}}(d)$ & $1.0750863337s$ & $1.3331823503s$ &
$0.1773717254s^2$\\ 
    \cline{2-5}
    Parallel Rhombuses & Sim & $1.0749140141s$ & $1.3318514164s$ &
$0.1764112787s^2$ \\ \hline
    Between two Long-Diag & $f_{D_{\rm LD}}(d)$ & $1.7570796617s$ & $3.3330145688s$
& $0.2456856311s^2$\\ 
    \cline{2-5}
     Adjacent Rhombuses & Sim & $1.7575618714s$ & $3.3319543736s$ &
$0.2429306419s^2$ \\ \hline
    Between two Short-Diag & $f_{D_{\rm SD}}(d)$ & $1.1150961004s$ &
$1.3332064091s$ & $0.0897670959s^2$ \\ 
    \cline{2-5}
     Adjacent Rhombuses & Sim & $1.1156689048s$ & $1.3341275970s$ &
$0.0894104918s^2$ \\ \hline
  \end{tabular}
  \label{tab:moment}
\end{table}

Table~\ref{tab:moment} lists the first two moments, and the variance of the random
distances in the four cases given in Section~\ref{sec:result}, and the
corresponding simulation results for verification purposes. 

\subsection{Polynomial Fits of Random Distances}

\begin{table}
  \caption{Coefficients of the Polynomial Fit and the Norm of Residuals (NR)}
  \centering
  \begin{tabular}{|c||c|c|}
    \hline
    PDF & Polynomial Coefficients & NR \\ \hline \hline
    & $10^8\times
\left[-0.000166~~0.002955~-0.024296~~0.122501~-0.423700\right.$ & \\  
    $f_{D_{\rm I}}(d)$ &
$1.065297~-2.013141~~2.916117~-3.273038~~2.858823~-1.941184$ & $0.095901$ \\ 
     & $1.018567~-0.408430~~0.123035~-0.027166~~0.004245~-0.000446$ & \\ 
     & $\left.0.000029~-0.00000116~~0.0000000892~-0.0000000000648\right]$ &
\\\hline
     & $10^5\times
\left[0.000019~-0.000513~~0.006194~-0.045754~~0.231329\right.$ & \\
     $f_{D_{\rm P}}(d)$ &
$-0.847940~~2.328313~-4.879792~~7.881444~-9.836044~~9.45481$ & $0.059485$ \\
     & $-6.93816~~3.828928~-1.553959~~0.4491195~ -0.088213~~0.010963$ & \\
     & $\left.-0.000782~~0.0000499~-0.0000000863~~-0.0000000025988\right]$ &
\\\hline
     & $10^4\times
\left[0.0000002~-0.0000074~~0.0001242~-0.0012723~~0.0089355\right.$ & \\ 
     $f_{D_{\rm LD}}(d)$ &
$-0.045576~~0.174538~-0.511656~~1.160176~-2.042559~~2.78728$ & $0.011340$ \\
     & $-2.927992~~2.339236~-1.395868~~0.606494~-0.185173~~0.037762$ & \\
     & $\left.-0.004740~~0.0003293~-0.0000098938~~0.0000000733372\right]$ & \\
\hline
     & $10^7\times
\left[-0.000022~~0.000493~-0.005116~~0.032591~-0.142436\right.$ & \\
     $f_{D_{\rm SD}}(d)$ &
$0.452238~-1.077816~~1.965337~-2.77016~~3.029748~-2.567377$ & $0.147836$ \\
     & $1.674771~ -0.831348~~0.308601~-0.083558~~0.015932~-0.002034$ & \\ 
     & $\left.0.000161~-0.00000698~~ 0.00000013178~-0.0000000006098\right]$ & \\
\hline
  \end{tabular}
  \label{tab:poly}
\end{table}

\begin{figure}
\centering
  \subfloat[Within a Single
Rhombus]{\includegraphics[width=0.5\columnwidth]{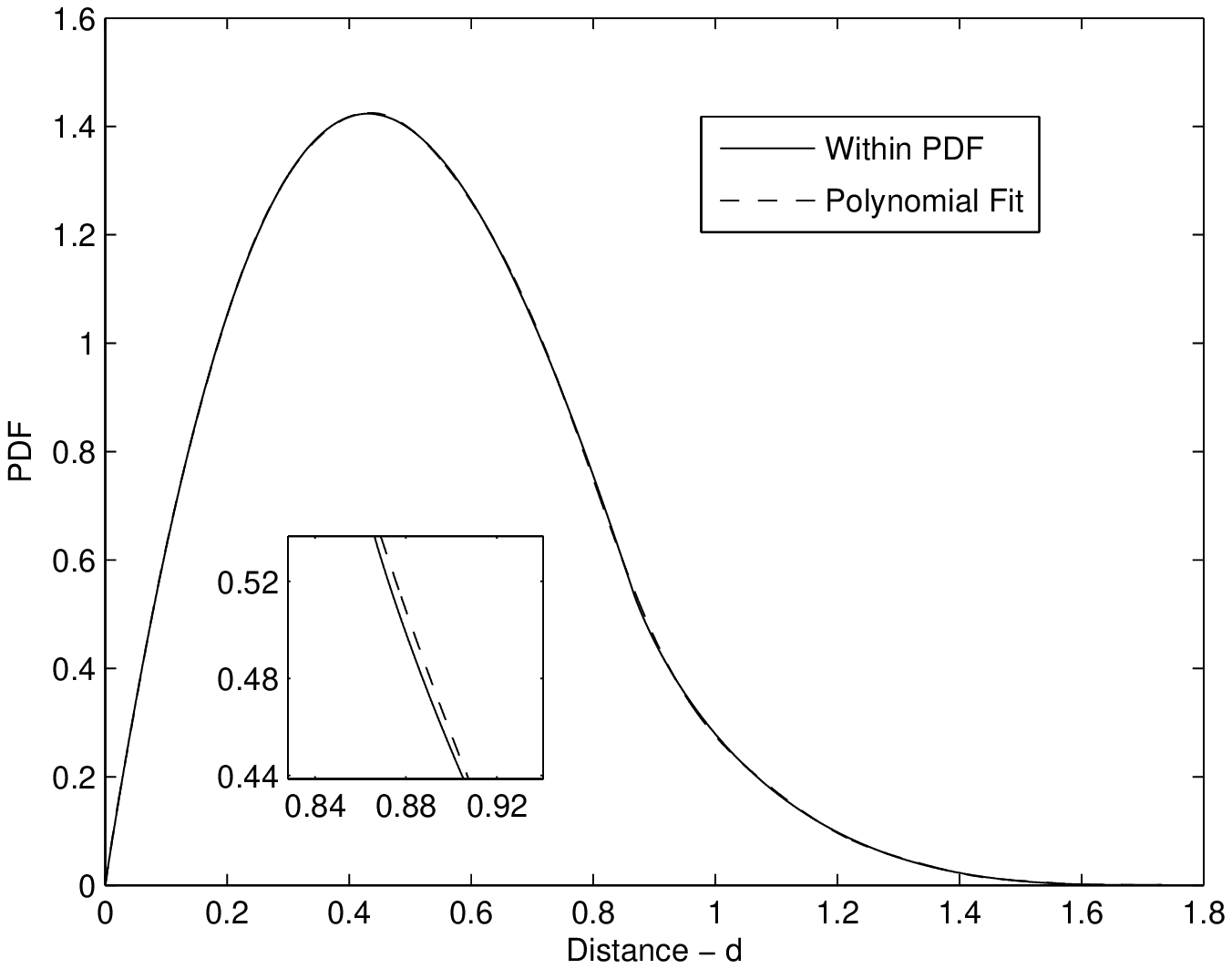}}
  \subfloat[Between two Parallel Adjacent
Rhombuses]{\includegraphics[width=0.5\columnwidth]{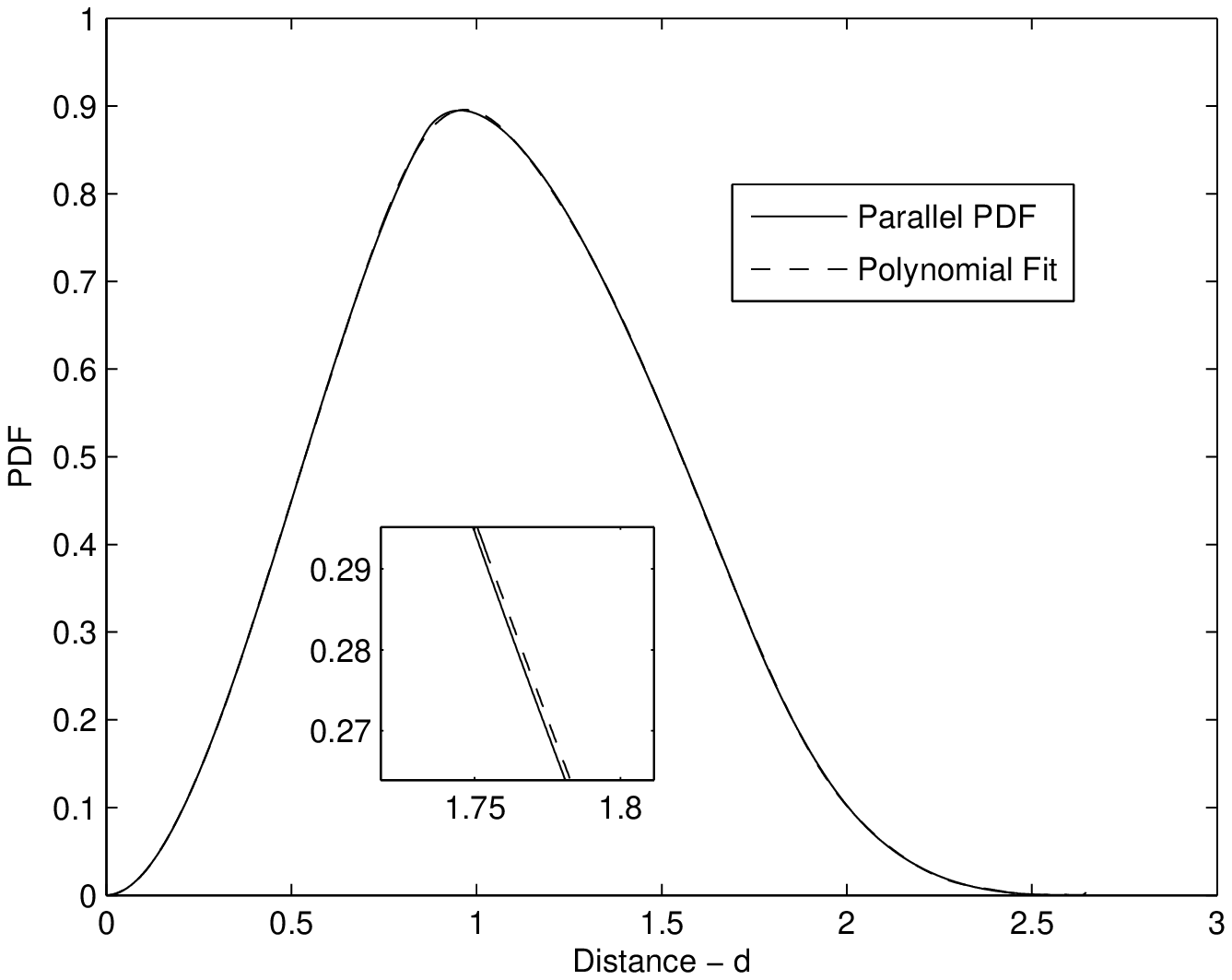}}\\
  \subfloat[Between two Long-Diag Adjacent
Rhombuses]{\includegraphics[width=0.5\columnwidth]{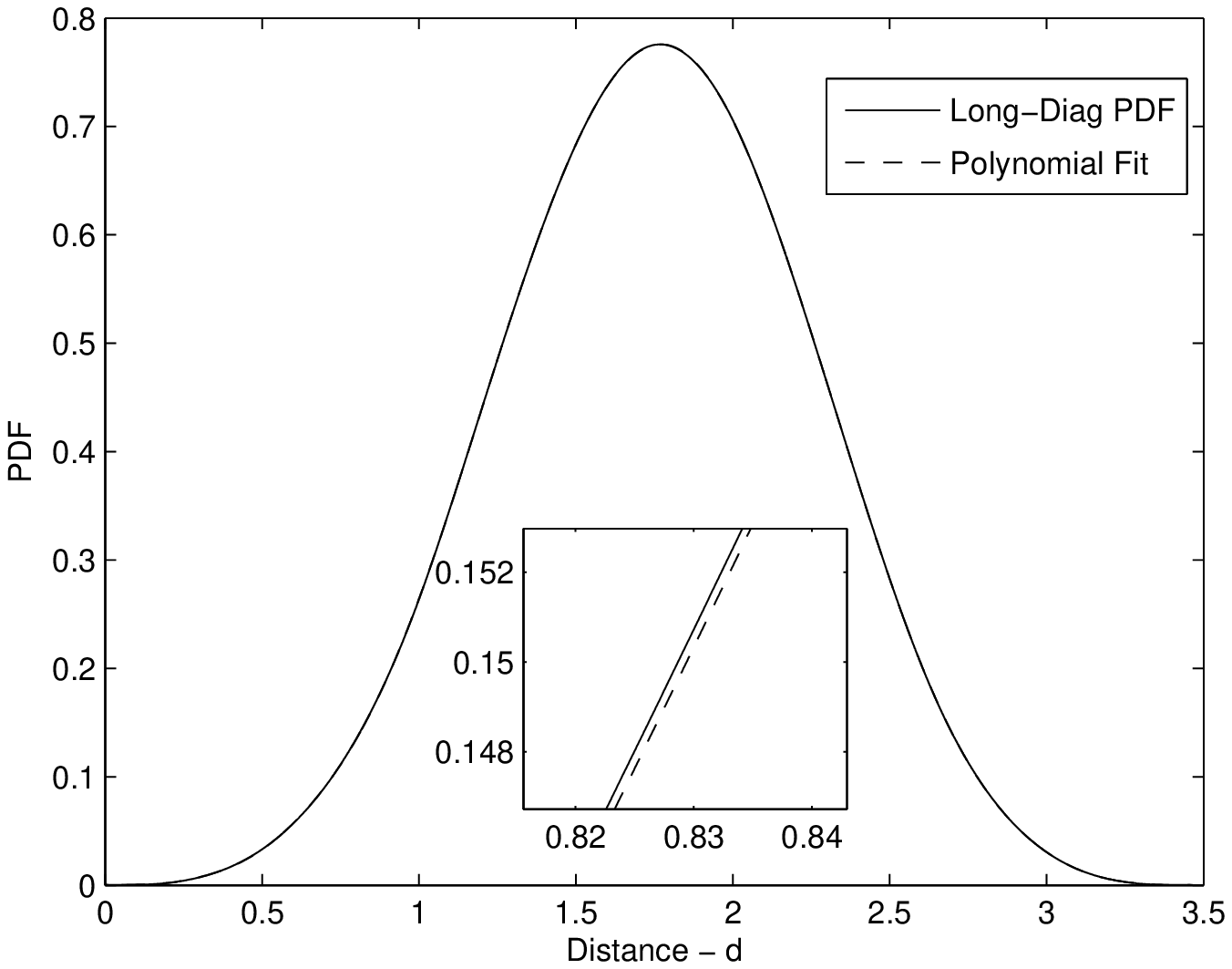}}
   \subfloat[Between two Short-Diag Adjacent
Rhombuses]{\includegraphics[width=0.5\columnwidth]{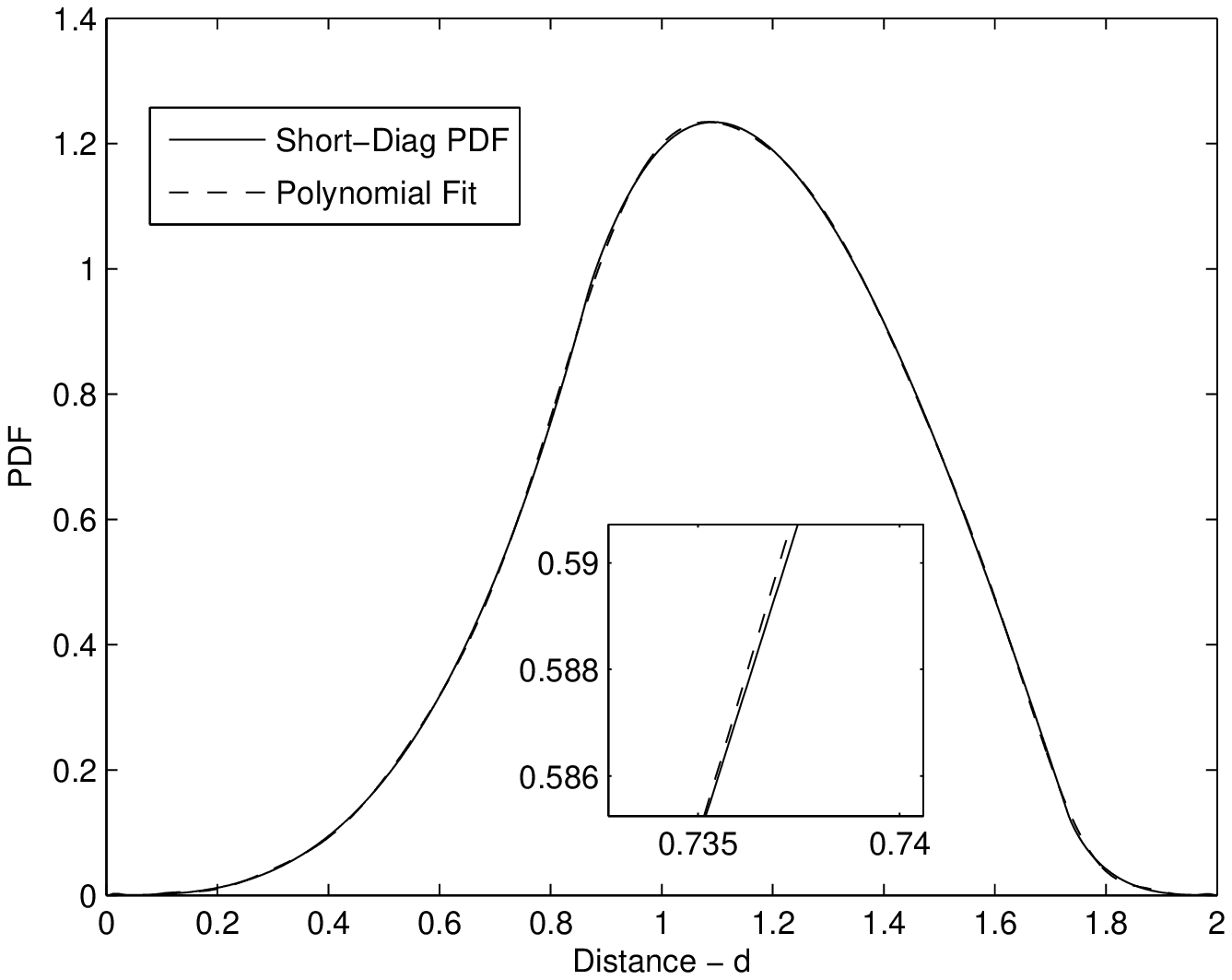}}
  \caption{Polynomial Fit of the Distance Distribution Functions Associated with Rhombuses.}
  \label{fig:rhombus_poly}
\end{figure}

Table~\ref{tab:poly} lists the coefficients of the degree-$20$ polynomial 
fits of the original PDFs given in Section~\ref{sec:result}, from $d^{20}$
 to $d^{0}$, and the corresponding norm of residuals. Figure~\ref{fig:rhombus_poly}
(a)--(d) plot the polynomials listed in Table~\ref{tab:poly} with the original PDFs. From the
figure, it can be seen that all the polynomials match closely with the original PDFs. These
high-order polynomials facilitate further manipulations of the distance
distribution functions, with a high accuracy.

\section{Conclusions}
\label{sec:conclude}

In this report, we gave the closed-form probability density functions
of the random distances associated with rhombuses. The correctness of the obtained
results has been validated by a recursion and a probabilistic sum, in
addition to simulation. The first two statistical moments, and the
polynomial fits of the density functions are also given for practical uses.

\section*{Acknowledgment}

This work is supported in part by the NSERC, CFI and BCKDF. The authors also
want to thank Dr. Lin Cai for the initial work involving squares, Dr. Aaron
Gulliver for posing the new problem associated with hexagons, and Dr. A.M.
Mathai, Dr. David Chu and Dr. T.S. Hale for their correspondence and
encouragement. The results presented here can lead to obtaining the distribution of random distances associated with hexagons.

\end{document}